\documentclass[10pt]{article}

\usepackage[utf8]{inputenc}
\usepackage[T2A]{fontenc}
\usepackage[english]{babel}

\usepackage{graphicx}
\usepackage{epsfig}

\usepackage[a4paper,
            bindingoffset=0.0in,
            left=0.9in,
            right=0.9in,
            top=0.9in,
            bottom=1.4in,
            footskip=.0in]{geometry}

\usepackage{lastpage}
\usepackage{fancyhdr} 
\usepackage{amssymb,amsthm,amsmath}
\usepackage{xcolor,paralist,hyperref,titlesec,fancyhdr,etoolbox}

\usepackage{epsf}
\usepackage{amsmath}
\usepackage{amsfonts}
\usepackage{amsthm}
\usepackage{amssymb}
\usepackage{enumerate}

\usepackage{pstricks}
\usepackage{lastpage}

\def\C{\mathbb{C}}
\def\R{\mathbb{R}}

\begin{document}
\begin{center}

{\Large

Application of Structured Matrices for Solving Hartree-Fock Equations

}

{\large\it

Ilgis Ibragimov

Elegant Mathematics LLC

ii@elegant-mathematics.com

}

\end{center}

\bigskip

{\Large\bf

\noindent
Abstract

}





\bigskip
\noindent
This work was originally published by the author in 1999 in a book \cite{1} and later became part of the author's doctoral thesis in 1999 \cite{2}. Since the original language of these works is not English, the author provides a translation of the key ideas of these publications in this work. In addition, the chapter related to numerical experiments was recalculated on modern computers and using contemporary benchmark datasets.

\bigskip

This article presents a novel approach to solving Hartree-Fock equations using Toeplitz and tensor matrices and bases based on regular finite elements.

The issues discussed include the choice of basis, the dependence of data volume and number of arithmetic operations on the number of basis functions, as well as the arithmetic complexity and accuracy of computing two- and four-center integrals.

The approach has been implemented in a software package, and results have been obtained that are in good agreement with theory.

\section{Introduction}






The difference between these eigenvalues characterizes the energies of electron emission or absorption by molecular systems. The eigenfunctions are responsible for the statistical properties of electrons, such as the probability of finding electrons in space around atomic nuclei at a given energy of the molecule.

The Hartree-Fock operator is defined in three-dimensional space and consists of the sum of three operators: $T$ -- the kinetic energy of electrons, $V_1$ -- the potential energy of electrons in the field of atomic nuclei, and $V_2$ -- the potential energy of electrons in the field created by all other electrons.

In terms of operator order, $T$ is larger than the other two operators, $V_1$ and $V_2$. Solving this problem for crystalline structures, where the boundary conditions are periodic in each Cartesian coordinate, is of particular interest. The discretization method based on Ritz \cite{19} is a common approach to solving this problem.

Let's introduce the following notations: $m$ -- the number of electron pairs, $N_n$ -- the number of atomic nuclei in the target system, and $N$ -- the size of the basis of test functions. The sought-after eigenfunctions of the Hartree-Fock operator are given by

\begin{equation}
\label{eq:3.4.2.1}
\psi_i=\sum_{j=1}^N c_{ji} \phi_j, \hspace*{10mm} i=1, \dots, m,
\end{equation}
where $c_{ji}$ are the unknown coefficients and $\phi_j$ are the test functions.
The form of the discrete Hartree-Fock operator is as follows:

\begin{equation}
\label{eq:3.4.2.2}
\sum_{p=1}^N (F_{pq} - \epsilon_k S_{pq}) c_{pk} =0,
\hspace*{10mm}
q=1, \dots, N;
k=1, \dots, m.
\end{equation}

\noindent where

$$
F_{pq} = <p|H|q> + \sum_{i,j=1; i \neq j}^N {P_{ij} {(<ij|pq> -
\frac{1}{2}<ip|jq>)}}
$$

\begin{equation}
\label{eq:3.4.2.2.1}
<p|H|q>=\int\limits \phi_p^*({\bf r}) \left\{ - \frac{1}{2} \left(
\frac{\partial^2}{\partial x^2}+\frac{\partial^2}{\partial y^2}+\frac{\partial^2}{\partial z^2} \right) +
\sum_{i=1}^{N_n} {\frac{{\bf Q}_i}{\left| {\bf R}_i - {\bf r} \right|}}
\right\} \phi_q({\bf r}) d{\bf r}
\end{equation}

\begin{equation}
\label{eq:3.4.2.2.2}
<ij|kl>=\int\limits
\frac{\phi_i^*({\bf r}_1) \phi_j({\bf r}_1) \phi_k^*({\bf r}_2)
\phi_l({\bf r}_2)} {|{\bf r}_1 - {\bf r}_2|} d{\bf r}_1 d{\bf r}_2
\end{equation}

$$
S_{ij}=<i|j>=\int\limits \phi_i^*({\bf r}) \phi_j({\bf r}) d{\bf r}
\hspace*{10mm}
P_{ij}=\sum_{k=1}^m {c^*_{ik} c_{jk}}.
$$




It is necessary to solve this system of nonlinear equations and find $c_{jk}$ and $\epsilon_k$.
In general, the matrices $F_{pq}$ and $S_{pq}$ are Hermitian.

Currently, there are numerous sets of test functions proposed \cite{19,11,12,20,24,26}. Broadly speaking, the types of bases can be divided into two categories: bases based on atomic orbitals and bases based on finite elements (FE bases).

Bases based on atomic orbitals are widely used, while FE bases have only recently been considered for this problem \cite{24,34}. This article specifically examines FE bases on a grid with identical cells. It has been found that sufficient accuracy for most calculations can be achieved with such bases.


\section{Bases Based on Regular Finite Elements}

We will solve the equation using the Ritz discretization method. Let $D$ be a domain in $\mathbb{R}^3$ in the form of a parallelepiped, such that all atomic nuclei are at a distance of at least $d$ from the complement of $D$.

In the domain $D$, we have a three-dimensional grid with a step size of $h$, such that along each Cartesian coordinate inside $D$, there are $n_i$, $n_j$, $n_k$ segments, respectively. This grid generates a set of identical cubes with a side length of $h$. Let each such cube have a three-dimensional index $i, j, k$.

Furthermore, let all atomic nuclei be located only at the centers of the cubes. Then, the coordinates of the nuclei can be expressed in terms of $h$ and the corresponding three-dimensional index of the cube in which the nucleus is located. We denote the coordinate vector as $R_{ijk}$.

Later, we will consider the case where the restriction on the placement of atomic nuclei needs to be lifted.

The basis functions ${ S_{ijk} }^{n_i, n_j, n_k}{i,j,k=1} \in \mathbb{R}^p$ are defined as $p$-smooth functions with finite support \cite{31}: $S{ijk}$ is nonzero in cubes with indices $p_i$, $p_j$, $p_k$ satisfying

$$\max ( |p_i - i|, |p_j - j|, |p_k - k|) \le p,$$
and in other cubes, $S_{ijk}$ is identically zero.

On the boundary, we will consider periodic conditions along each Cartesian coordinate, although all subsequent derivations can be easily generalized to zero boundary conditions. Also, let all $S_{ijk}$ be obtainable from any other $S_{i'j'k'}$ through parallel translation.

In most cases, two additional conditions will be considered:
\begin{itemize}
\item [1)] Orthogonality of all $S_{ijk}$.
\item [2)] Symmetry of all $S_{ijk}$ with respect to three planes passing through the center of the cube with index $i, j, k$, where each of these planes is parallel to one of the faces of the domain $D$.
\end{itemize}



The discretization in this basis is possible due to the Hermitian nature of the original operator \cite{38}. The details associated with this will not be discussed in this work.

Since the basis consists of functions with finite support, some integrals in formula (\ref{eq:3.4.2.2}) will be identically zero. Thus, the integrals

\begin{eqnarray}
\label{eq:4.0.1}
\int\limits_D S^*_{ijk}(\bar r) S_{i'j'k'}(\bar r) d \bar r,
\end{eqnarray}

\begin{eqnarray}
\label{eq:4.0.2}
- \frac{1}{2} \int\limits_D S^*_{ijk}(\bar r)
  \left[ \frac{\partial^2}{\partial x^2} + \frac{\partial^2}{\partial y^2} + \frac{\partial^2}{\partial z^2} \right]
  S_{i'j'k'}(\bar r) d \bar r,
\end{eqnarray}

\begin{eqnarray}
\label{eq:4.0.3}
\int\limits_D S^*_{ijk}(\bar r)
  \left[ \sum_{q=1}^{N_n} {\frac{{\bf Q}_q}{\left| \bar R_q - \bar r \right|}}
  \right] S_{i'j'k'}(\bar r) d \bar r
\end{eqnarray}
%
%
%
%
%
are identically zero in the case of zero boundary conditions, if $\max ( |i - i'|, |j - j'|, |k - k'|) > p$, and in the case of periodic boundary conditions along each Cartesian coordinate, if

$$\max ( {\rm MOD}(i - i'+n_i, n_i), {\rm MOD}(j - j'+n_j, n_j), {\rm MOD}(k - k'+n_k, n_k)) > p.$$
Here, ${\rm MOD}(a,b)$ denotes the operation of taking the remainder of dividing $a$ by $b$.

For the sake of simplicity, we will consider only periodic boundary conditions, as well as both additional conditions: orthogonality and symmetry of all $S_{ijk}$. Then, (\ref{eq:4.0.1}) is equal to $1$ when $i=i'$, $j=j'$, $k=k'$, and it is equal to $0$ in all other cases. For integrals of the form (\ref{eq:4.0.2}) and (\ref{eq:4.0.3}), the equalities

\begin{eqnarray}
\label{eq:4.0.4}
\nonumber
\int\limits_D S^*_{ijk}(\bar r)
  \left[ \frac{\partial^2}{\partial x^2} + \frac{\partial^2}{\partial y^2} + \frac{\partial^2}{\partial z^2} \right]
  S_{{\rm MOD}(i+i',n_i), {\rm MOD}(j+j',n_j), {\rm MOD}(k+k',n_k)}(\bar r) d \bar r=
\\
\int\limits_D S^*_{i''j''k''}(\bar r)
  \left[ \frac{\partial^2}{\partial x^2} + \frac{\partial^2}{\partial y^2} + \frac{\partial^2}{\partial z^2} \right]
  S_{{\rm MOD}(i''+i',n_i), {\rm MOD}(j''+j',n_j), {\rm MOD}(k''+k',n_k)}(\bar r) d \bar r,
\end{eqnarray}

\begin{eqnarray}
\label{eq:4.0.5}
\nonumber
\int\limits_D S^*_{ijk}(\bar r)
  \left[ {\frac{1}{\left| \bar R_{ijk} - \bar r \right|}} \right]
  S_{{\rm MOD}(i+i',n_i), {\rm MOD}(j+j',n_j), {\rm MOD}(k+k',n_k)}(\bar r) d \bar r=
\\
\int\limits_D S^*_{i''j''k''}(\bar r)
  \left[ {\frac{1}{\left| \bar R_{i''j''k''} - \bar r \right|}} \right]
  S_{{\rm MOD}(i''+i',n_i), {\rm MOD}(j''+j',n_j), {\rm MOD}(k''+k',n_k)}(\bar r) d \bar r,
\end{eqnarray}
%

%
hold for any $i, i', i''$, $j, j', j''$, $k, k', k''$. It can be observed that the number of different values of integrals of the form (\ref{eq:4.0.4}) is of the order of $p^3$, and the number of integrals of the form (\ref{eq:4.0.5}) is of the order of $n_i n_j n_k$.

Computing integrals of the form (\ref{eq:3.4.2.2.2}) will be more complex.

\begin{eqnarray}
\label{eq:4.0.6}
\int\limits_D \int\limits_D \frac{S^*_{i_1j_1k_1}(\bar r) S_{i_2j_2k_2}(\bar r) 
                    S^*_{i_3j_3k_3}(\bar r') S_{i_4j_4k_4}(\bar r')}
                   {\bar r - \bar r'} d \bar r d \bar r'.
\end{eqnarray}


Due to the basis being constructed with functions of finite support, the integral (\ref{eq:4.0.6}) is identically zero either when

\begin{eqnarray}
\label{eq:4.0.7}
\max ( {\rm MOD}(i_1 - i_2 + n_i, n_i), {\rm MOD}(j_1 - j_2 + n_j, n_j), {\rm MOD}(k_1 - k_2 + n_k, n_k)) > p,
\end{eqnarray}
or when
\begin{eqnarray}
\label{eq:4.0.8}
\max ( {\rm MOD}(i_3 - i_4 + n_i, n_i), {\rm MOD}(j_3 - j_4 + n_j, n_j), {\rm MOD}(k_3 - k_4 + n_k, n_k)) > p.
\end{eqnarray}

Since the basis functions are symmetric, the equality

$$
\int\limits_D \int\limits_D \frac{S^*_{i_1j_1k_1}(\bar r) S_{i_2j_2k_2}(\bar r) 
                    S^*_{i_3j_3k_3}(\bar r') S_{i_4j_4k_4}(\bar r')}
                   {\bar r - \bar r'} d \bar r d \bar r'=
\int\limits_D \int\limits_D \frac{S^*_{i_3j_3k_3}(\bar r) S_{i_4j_4k_4}(\bar r) 
                    S^*_{i_1j_1k_1}(\bar r') S_{i_2j_2k_2}(\bar r')}
                   {\bar r - \bar r'} d \bar r d \bar r'.
$$
%
holds. Finally, since any basis function can be obtained from another through parallel translation, the equality

$$
\int\limits_D \int\limits_D \frac{S^*_{i_1j_1k_1}(\bar r) S_{i_2j_2k_2}(\bar r) 
                    S^*_{i_3j_3k_3}(\bar r') S_{i_4j_4k_4}(\bar r')}
                   {\bar r - \bar r'} d \bar r d \bar r'=
\int\limits_D \int\limits_D \frac{S^*_{i_5j_5k_5}(\bar r) S_{i_6j_6k_6}(\bar r) 
                    S^*_{i_7j_7k_7}(\bar r') S_{i_8j_8k_8}(\bar r')}
                   {\bar r - \bar r'} d \bar r d \bar r',
$$
holds under the conditions

$$\begin{tabular}{ll}
${\rm MOD}(i_1+i_3, n_i)={\rm MOD}(i_5+i_7, n_i)$, & ${\rm MOD}(i_2+i_4, n_i)={\rm MOD}(i_6+i_8, n_i)$, \\
${\rm MOD}(j_1+j_3, n_j)={\rm MOD}(j_5+j_7, n_j)$, & ${\rm MOD}(j_2+j_4, n_j)={\rm MOD}(j_6+j_8, n_j)$, \\
${\rm MOD}(k_1+k_3, n_k)={\rm MOD}(k_5+k_7, n_k)$, & ${\rm MOD}(k_2+k_4, n_k)={\rm MOD}(k_6+k_8, n_k)$.
\end{tabular}$$

Let us denote
$$
t(i_2, j_2, k_2, i_3, j_3, k_3, i_4, j_4, k_4)=
\int\limits_D \int\limits_D \frac{S^*_{0,0,0}(\bar r) S_{i_2j_2k_2}(\bar r) 
                    S^*_{i_3j_3k_3}(\bar r') S_{i_4j_4k_4}(\bar r')}
                   {\bar r - \bar r'} d \bar r d \bar r'.
$$





Considering (\ref{eq:4.0.7}) and (\ref{eq:4.0.8}), it is easy to notice that there are a number of different values of $t$ on the order of $p^6 n_i n_j n_k$. The coefficient $p^6$ is greatly overestimated -- due to symmetry, most of these integrals are equal to each other. In most cases, we will limit ourselves to the values $p=3$, and for this case, the constant in front of $n_i n_j n_k$ will be equal to $28$.

We will perform a change of variables in all integrals (\ref{eq:4.0.2}), (\ref{eq:4.0.3}), and (\ref{eq:4.0.6}) such that in the new coordinates, the grid spacing is identically equal to one. Then, in the new coordinates, the coefficient $\frac{1}{h}$ will appear in front of (\ref{eq:4.0.3}) and (\ref{eq:4.0.6}), while the coefficient $\frac{1}{h^2}$ will appear in front of (\ref{eq:4.0.2}).

It can also be shown that if all first and second partial derivatives of $S_{ijk}$ are bounded, then all integrals in the new coordinate system are bounded by a constant that does not depend on the problem size.

To summarize the results obtained, for solving the discretized Hartree-Fock equations using the Ritz method, it is necessary to compute three types of integrals: (\ref{eq:4.0.2}), (\ref{eq:4.0.3}), and (\ref{eq:4.0.6}). The number of different integrals of the form (\ref{eq:4.0.2}) does not depend on the problem size. Moreover, for piecewise polynomial basis functions, such integrals can be analytically computed without much difficulty. The number of different integrals of the form (\ref{eq:4.0.3}) and (\ref{eq:4.0.6}) is linearly dependent on the problem size.


\subsection{Integration Computation}


Let us consider the computation of integrals of the form

\begin{eqnarray}
\label{eq:4.1.2}
\int\limits_{\Omega} \frac{s(\bar r)}{||\bar r - \bar R||_2} d \bar r
\end{eqnarray}
and
\begin{eqnarray}
\label{eq:4.1.3}
\int\limits_{\Omega} \int\limits_{\Omega}
  \frac{s_1(\bar r_1) s_1(\bar r_2)}{||\bar r_1 - \bar r_2 + \bar t||_2}
  d \bar r_1 d \bar r_2,
\end{eqnarray}
%
%
where $\bar r$, $\bar r_1$, $\bar r_2$, $\bar t$, $\bar R$ are three-dimensional vectors, $\Omega$ is a parallelepiped-shaped region, $s$, $s_1$, $s_2$ are functions that are tensor products of functions with a single argument: $s(\bar x) = f_1(x_1) f_2(x_2) f_3(x_3).$


Since each $s(\bar x)$ is a product of two basis functions, whose domain is a cube with side length $h$, the integration domain $\Omega$ cannot be larger than a cube with side length $h$. The coordinates of the vectors $\bar R$ and $\bar t$ take integer values from $0$ to $n-1$, where $n$ is the size of the problem.

 
Note that under the conditions $R_x \leq h$, $R_y \leq h$, $R_z \leq h$, the integrand in (\ref{eq:4.1.2}) can take infinitely large absolute values.
At the same time, when $||\bar R||_2 >> h$, the integrand is bounded, varies little, and the value of the integral (\ref{eq:4.1.2}) tends to
$\displaystyle \frac{1}{||\bar R||_2} \int\limits_{\Omega} s(\bar x) d \bar x$ при $h \to \infty$.


Similar reasoning can be applied to the integral (\ref{eq:4.1.3}), using $\bar t$ instead of $\bar R$.


Let us also mention a few properties of the functions $s$:
\begin{itemize}
\item [1)] $s$ is a piecewise polynomial function;
\item [2)] since $s$ is the product of two basis functions, and by construction, the basis functions are symmetric with respect to $x_i=0$, then $s$ is symmetric with respect to $x_1=0$, $x_2=0$, and $x_3=0$;
\item [3)] $s$ vanishes on the boundary of the integration domain $\Omega$, except when $s$ is piecewise constant.
\end{itemize}



It is necessary to compute the integrals (\ref{eq:4.1.2}) and (\ref{eq:4.1.3}) with an absolute accuracy of at least $\frac{\eta}{N}$, where $\eta$ is the absolute accuracy of eigenvalue computation.

Next, we will discuss several techniques that allow computing such integrals with the desired accuracy using a reasonable number of arithmetic operations.

\bigskip


\noindent
{\bf Theorem 1.} For any $p$ and $q$, which are tensor products of functions of a single variable, and for a parallelepiped-shaped region $\Omega$, there exist $s$ and $\Omega'$ such that
\begin{eqnarray}
\label{eq:4.1.4}
\int\limits_{\Omega} \int\limits_{\Omega}
  \frac{p(\bar x) q(\bar y)}{||\bar x - \bar y + \bar t||_2}
  d \bar x d \bar y =
\int\limits_{\Omega'} \frac{z(\bar r)}{||\bar r - \bar t||_2} d \bar r.
\end{eqnarray}


\noindent
{\bf Proof.} The integral $\displaystyle\int\limits_0^h f_y(y) dy \int\limits_0^h \omega(x-y) f_x(x) dx$ upon the variable substitutions $u=x-y$ and $v=x+y$ transforms into
$$
\int\limits_{-h}^0 \omega(u)
  \int\limits_{h+u}^{h-u} f_x\left(\frac{v+u}{2}\right)f_y\left(\frac{v-u}{2}\right) dv du +
\int\limits_0^h \omega(u)
  \int\limits_u^{2h-u} f_x\left(\frac{v+u}{2}\right)f_y\left(\frac{v-u}{2}\right) dv du.
$$

Let us denote
$$
g(u) = \left\{ \begin{array}{lcl}
u \le 0 & : & \displaystyle \int\limits_{h+u}^{h-u} f_x\left(\frac{v+u}{2}\right)f_y\left(\frac{v-u}{2}\right) dv \\
\\
u >   0 & : & \displaystyle \int\limits_u^{2h-u} f_x\left(\frac{v+u}{2}\right)f_y\left(\frac{v-u}{2}\right) dv \\
\end{array} \right.
$$
%
Then
\begin{eqnarray}
\label{eq:4.1.5}
\int\limits_0^h f_y(y) dy \int\limits_0^h \omega(x-y) f_x(x) dx =
\int\limits_{-h}^h \omega(u) g(u) du.
\end{eqnarray}
%
%
Let the region $\omega$ be given as a parallelepiped with sides $h_1 \times h_2 \times h_3$. According to the conditions
$p(\bar x) = p_1(x_1) p_2(x_2) p_3(x_3)$ and $q(\bar y) = q_1(y_1) q_2(y_2) q_3(y_3)$, and using (\ref{eq:4.1.5}),
we obtain (\ref{eq:4.1.4}), where $\Omega'$ is a parallelepiped with sides $2 h_1 \times 2 h_2 \times 2 h_3$.
This completes the proof.

\bigskip



In most cases, the function $z$ in (\ref{eq:4.1.4}) can be constructed analytically: if $p$ and $q$ are basis splines of degree $k$, then $z$ is also a basis spline of degree $2k+1$.

Thus, the set of problems has significantly simplified, and it is necessary to compute only integrals of the form (\ref{eq:4.1.2}), as the basis functions
can be chosen in such a way that the indefinite integrals of the form $\displaystyle \int\limits \int\limits \int\limits s_1(\bar x + \bar y) s_2(\bar x + \bar y) d \bar x$
can be computed analytically.

\bigskip


\noindent
{\bf Lemma 1.} The indefinite integral $\displaystyle \int\limits \int\limits \int\limits \frac{x_1^{n_1} x_2^{n_2} x_3^{n_3}}{R} dx_1 dx_2 dx_3,$ where
$\displaystyle R=\sqrt{x_1^2 + x_2^2 + x_3^2}$, and $n_1$, $n_2$, $n_3$ are positive integers, can be represented as
$$
p_1 R +
p_2 \ln (x_1 + R) +
p_3 \ln (x_2 + R) +
p_4 \ln (x_3 + R) +
p_5 {\rm arctg} \left(\frac{x_1 x_2}{x_3 R}\right) +
p_6 {\rm arctg} \left(\frac{x_2 x_3}{x_1 R}\right) +
p_7 {\rm arctg} \left(\frac{x_1 x_3}{x_2 R}\right),
$$
%
%
where $p_1, \dots, p_7$ are polynomials in $x_1$, $x_2$, and $x_3$ of degree not exceeding $3 \max(n_1, n_2, n_3) + 2$.


The proof of the lemma is based on direct integration and is not presented here due to its complexity.

\bigskip




The lemma supports the idea of analytical integration when the basis functions are piecewise polynomials. However, as $R \to \infty$, the definite integral approaches $\frac{1}{R}$, while the indefinite integral has a magnitude of $R^{3k+2}$, where $k$ is the maximum degree of the polynomials in the integrand. Therefore, the absolute integration error will be no less than $\eta R^{3(k+1)}$. For $k \geq 2$, such an approach becomes practically unacceptable for computations due to the large error.

For $k = 0, 1$, analytical computation of integrals remains stable for small $R$, and for $R < 10$, it yields an error no greater than $10^{-9}$. Calculating a single integral of this type requires about 200 arithmetic operations and 56 evaluations of irrational functions such as $\ln$, ${\rm arctg}$, $\sqrt{.}$

Unfortunately, for $k=2$, the analytical formulas already involve polynomials up to degree 8, and the number of such polynomials is very large, around 150. The absolute accuracy of integration obtained experimentally (using numerical integration with Gaussian quadrature as the ``exact'' value) reached only $10^{-3}$ for $R=0$, indicating that the accuracy for $R>0$ can only be worse, which is completely unacceptable for our purposes.



In most cases, the function $s(\bar r)$ in (\ref{eq:4.1.2}) does not belong to the class of piecewise linear functions, so it was necessary to find another way to compute such integrals.

Using Gaussian quadrature formulas without step selection requires enormous computational effort in the vicinity of the point $||\bar R||_2=0$. Gaussian quadrature formulas with step selection for three-dimensional problems do not significantly reduce the number of arithmetic operations. Therefore, it is suggested to approximate $s(\bar r)$ using piecewise linear functions. In the obtained partition, each parallelepiped contains an analytically computable integral. The accuracy of such integration is the product of the approximation accuracy of $s(\bar r)$ with piecewise linear functions $l(\bar r)$ and the maximum value of the integral $\displaystyle \int\limits_{\Omega} \frac{d \bar x}{||\bar R - \bar x||_2} $ over all possible $\bar R$. It is obvious that the maximum is reached when $R$ is the center of the parallelepiped $\Omega$.


Thus, we obtain the final estimate:
$$
|s(\bar r) - f(\bar r)| \int\limits_{\Omega} \frac{d \bar x}{||\bar R - \bar x||_2}.
$$


The value $\displaystyle \int\limits_{\Omega} \frac{d \bar x}{||\bar R - \bar x||_2}$ is analytically computable and is bounded above by the number $3$ for $\Omega$ in the form of a cube with sides $1 \times 1 \times 1$.


Non-uniform partitioning can be used to reduce $\max |s(\bar r) - f(\bar r)|.$


This problem is straightforward since $f_i$ and $s_i$ are piecewise polynomial functions. Moreover, if $s_i(-h)=0$ and $s_i(h)=0$, then
\begin{eqnarray}
\label{eq:4.1.7}
\int\limits_{-h}^h \int\limits_{-h}^h \int\limits_{-h}^h \frac{s(\bar x)}{||\bar R - \bar x||_2}
dx_1 dx_2 dx_3 =
\sum_{i,j,k=1}^n F_{ijk} p_{ijk},
\end{eqnarray}
%
%
where $F_{ijk}$ depends only on $R$ and not on $s$, and $p_{ijk}$ depends only on the form of $s$ and not on $R$. This is achieved because $F_{ijk}$ are quantities of integrals that are computed analytically, and $p_{ijk}$ are their weights, which depend only on $s$.



If multiple integrals need to be computed in the same domain, using formula (\ref{eq:4.1.7}) significantly simplifies the computation since the most computationally intensive part, computing $F_{ijk}$, is performed once for multiple integrals.

Thus, we have shown how to compute integrals when the value of $||\bar R||_2$ is in the vicinity of zero. Now, we need to understand how to compute integrals for which $||\bar R||_2 >>0$.


As mentioned above, for such integrals, the subintegral function changes only slightly. Therefore, it was proposed to expand the function $\displaystyle \frac{1}{||\bar R - \bar x||_2} $ in a Taylor series with respect to all Cartesian coordinates $\bar x$ at the point $\bar R$. In this case, a series in negative powers of the function $||\bar R||_2$ is obtained, with coefficients $p_i, i=1, \dots$ in the form of polynomials of the coordinates $\bar R$ and $\bar x$. The first coefficient is equal to one. The estimate of the remainder term for the expansion up to the $k$-th degree corresponds to the expression
$$
{\tt const} \int\limits_{-h}^h \int\limits_{-h}^h \int\limits_{-h}^h \frac{dx_1 dx_2 dx_3}
                                                {||\bar R||_2^{k+1} (k+1)!},
$$
%
%


An interesting fact was also discovered: in this expansion, the coefficients of even powers in $||\bar R||_2$ only contain odd powers of the coordinates $\bar x$.


If integrals of the form (\ref{eq:4.1.2}) with $||\bar R||_2 >> 0$ are computed using such a Taylor expansion, then (\ref{eq:4.1.2}) can be rewritten as:
$$
\int\limits_{\Omega} \frac{s(\bar r)                   }{||\bar R||_2} \bar r +
\int\limits_{\Omega} \frac{s(\bar r) p_2(\bar r,\bar R)}{||\bar R||_2^2} \bar r +
\int\limits_{\Omega} \frac{s(\bar r) p_3(\bar r,\bar R)}{||\bar R||_2^3} \bar r + \dots =
$$
\begin{eqnarray}
\label{eq:4.1.9}
\frac{1}{||\bar R||_2}   \int\limits_{\Omega} s(\bar r)                    d \bar r +
\frac{1}{||\bar R||_2^2} \int\limits_{\Omega} s(\bar r) p_2(\bar r,\bar R) d \bar r +
\frac{1}{||\bar R||_2^3} \int\limits_{\Omega} s(\bar r) p_3(\bar r,\bar R) d \bar r + \dots
\end{eqnarray}
%
As can be seen, the coefficients of the expansion in powers of $||\bar R||_2$ are integrals of products of basis functions and certain polynomial functions and, in most cases, can be computed analytically. Therefore, all these coefficients can be precomputed.


Furthermore, since the functions $s(\bar r)$ are symmetric by assumption, all integrals involving even powers of $||\bar R||_2$ evaluate to zero. The accuracy estimate of such integration is:
$$
\frac{1}{||\bar R||_2^k} \int\limits_{\Omega} s(\bar r) p_k(\bar r) d \bar r.
$$

%
%

Thus, an interesting problem arises: there are two alternative algorithms for computing integrals that are determined by the parameter $R$ and the form of the basis functions. Each of these algorithms depends on an additional parameter: for the first algorithm, it is $k_1$ --- the number of nodes in the partition, and for the second algorithm, it is $k_2$ --- the maximum degree of the Taylor expansion. For each of these algorithms, there exist:
\begin{itemize}
\item a priori functions $\epsilon_1(k_1, R)$ and $\epsilon_2(k_2, R)$ that can be used to determine the accuracy of computing such an integral for the corresponding algorithm;
\item a priori functions $f_1(k_1, R)$ and $f_2(k_2, R)$ that can be used to determine the computational costs required to compute this integral with the accuracy of $\epsilon_1(k_1, R)$ and $\epsilon_2(k_2, R)$, respectively.
\end{itemize}
It is necessary to correctly choose the algorithm and the corresponding number $k_1$ or $k_2$ for computing integrals with an accuracy of at least $\eta$.



Since all these functions are monotonic with respect to $||\bar R||_2$ over the entire domain, the solution to this problem will be the condition: ``there exists $\alpha$ such that when $||\bar R||_2 < \alpha$, the first algorithm with a constant value of $k_1$ is used. Otherwise, the second algorithm with certain values of $k_2$ is used.''

As practice has shown, for $||\bar R||_2 > \alpha$, for the majority of integrals, the accuracy obtained with $k_2 \le 3$ is sufficient. That is, in the expansion (\ref{eq:4.1.9}), only the first one or two nonzero terms are used. In this case, the computational costs for calculating a single integral correspond to a constant of the order of 20 arithmetic operations.

\subsection{Properties of the Discretized Hartree-Fock Functional}




Currently, there are numerous algorithms available for minimizing multivariable functions \cite{32,33}. Therefore, let us first consider how to compute the value of a function at a point, as well as its first and second partial derivatives with respect to all the variables. Only then will we discuss the minimization algorithm. We will also focus on approximate computations of the aforementioned operations, assuming that the computational costs and/or the amount of data for temporary and permanent storage are significantly smaller than those required for exact computations. This may be necessary when constructing a preconditioner.

Let $N=n_i n_j n_k$ be the size of the problem, and $m$ be, as before, the number of electron pairs. Each eigenfunction corresponding to an electron pair is represented by a vector $c_i \in \mathbb{C}^N$, where $i=1, \dots, m$.

The number of integrals required to solve the Hartree-Fock equations is of the order of $N$. The number of unknowns in the problem is $Nm$ complex numbers.

\subsection{Computation of the Value of the Discretized Hartree-Fock Energy Functional}


Let us express the Hartree-Fock energy in matrix form as
$$E = T_e + V_{en} + V_{ee},$$
and consider each term separately.

%





\subsubsection{Kinetic Energy of All Electrons}

The kinetic energy of all electrons is given by:
$$T_e = \frac{1}{h^2} \sum_{s=1}^m c_s^* A c_s.$$
Here, $A$ is defined by a parameter set of order $p^3$, which is independent of $N$. Furthermore, such parameters do not depend on the size of the matrix but only on the type of basis functions.

$A$ is a three-level positive-definite band matrix. Each level of the matrix has the same elements on its diagonals, forming circulants. The number of such diagonals is $(2p+1)^3$. Additionally, the matrix is symmetric.

It is known that an arbitrary circulant matrix (including multi-level circulants) can be represented as a product: $F^* D F$, where $F$ is the Fourier matrix and $D$ is a diagonal matrix. Multiplying a vector of size $N$ by the Fourier matrix requires only about $N\log_2 N$ arithmetic operations. Therefore, computing the kinetic energy of all electrons requires approximately $m N \min(\log_2 N, (2p+1)^3)$ arithmetic operations.

For small values of $p$, it is more efficient to use explicit matrix-vector multiplication, while for large values, the Fast Fourier Transform (FFT) is preferred.

\subsubsection{Potential Energy of All Electrons in the Field of All Atomic Nuclei}


The potential energy of all electrons in the field of all atomic nuclei is given by:
$$V_{en} = - \frac{2}{h} \sum_{s=1}^m c_s^* B c_s.$$
%
%
Here, $B$ is a three-level positive-definite band matrix with the number of nonzero diagonals equal to $(2p+1)^3$.


Furthermore, the matrix $B$ can be represented as a sum:
\begin{eqnarray}
\label{eq:4.2.1.1}
B=\sum_{s=1}^{N_n} Q_s P(i_s, j_s, k_s) W P(n_i-i_s, n_j-j_s, n_k-k_s),
\end{eqnarray}
%
%
where $W$ is a three-level positive-definite band matrix with the number of nonzero diagonals equal to $(2p+1)^3$, consisting of the integrals (\ref{eq:4.1.2}), and the matrices $P(i, j, k)$ are permutation matrices that shift the element at position $(0,0,0)$ to position $(i,j,k)$. The indices $i_s, j_s, k_s$ correspond to the cube number containing the $s$-th atomic nucleus.



The choice of using the representation (\ref{eq:4.2.1.1}) instead of the matrix $B$ depends on the size of the problem. For large $N$, a significant portion of the off-diagonal elements of the matrix $W$ tend to zero, allowing us to neglect them. In this case, the total memory required to store such a compact representation and the number of arithmetic operations for multiplying a vector by such a matrix can be significantly smaller than when using the matrix $B$. The correct choice can only be made by knowing $N$ and $p$.

For large $N$, the asymptotic number of operations required to compute $V_{en}$ is on the order of $N m$.

\subsubsection{Potential Energy of All Electrons in the Field Generated by the Electrons Themselves}

The potential energy of all electrons in the field generated by the electrons themselves is given by:

\begin{equation}
\label{eq:4.2.1.2}
V_{ee} = \frac{1}{h} \sum_{s,q=1}^m \left\{
  2 c_s^* H^*(c_s) T H(c_q) c_q - c_s^* H^*(c_q) T H(c_s) c_q \right\},
\end{equation}
%
%
Here, $T$ is a three-level block-symmetric circulant matrix with nonnegative elements, constructed from the values of integrals of the form (\ref{eq:4.1.3}). $H(c)$ is a three-level block-band matrix whose elements are formed from the vector $c$, where the $(i,j,k)$-th block represents a vector of size $p^3$ composed of the numbers $c_{i+i', j+j', k+k'}$, where $i',j',k' \in [0, p]$.

When using the Fast Fourier Transform (FFT), the computational cost for computing $V_{ee}$ is on the order of ${\tt const} m^2 N \log_2 N$, where ${\tt const} \simeq p^6$.

\subsubsection{Computation of First Derivatives of the Discretized Hartree-Fock Energy Functional}


Computing the first derivatives of the discretized Hartree-Fock energy functional is algorithmically not much different from computing the functional itself. Therefore, only the final formula is provided here (using the notations from the previous section):
\begin{eqnarray}
\label{eq:4.2.2.1}
dE_i = \frac{2}{h^2} A c_i - \frac{4}{h} B c_i +
\frac{4}{h} \sum_{s=1}^m \left\{
  2 H^*(c_i) T H(c_s) c_s - H^*(c_s) T H(c_s) c_i \right\} -
\epsilon_i c_i,
\end{eqnarray}
%
%
Here, $dE_i$, $i=1, \dots, m$, are the vectors of first derivatives of the functional with respect to the coordinates of the $i$-th vector.


The number of arithmetic operations required to compute all such vectors is no more than that needed for computing the value of the discretized Hartree-Fock energy functional itself.

\subsubsection{Computation of Second Derivatives of the Discretized Hartree-Fock Energy Functional}


Due to the cumbersome expression, here we present the form of the matrix of second derivatives of the discretized Hartree-Fock energy functional only for the case of $p=0$:
\begin{eqnarray*}
\nonumber
d^2 E_{ii} = \frac{2}{h^2} A - \frac{4}{h} B +
\frac{4}{h} \left\{ 3 {\rm diag} (c^*_i)T {\rm diag} (c_i)- {\rm diag} (T {\rm diag} (c^*_i)c_i) \right\} +
\\
\frac{4}{h} \sum_{s=1}^m \left\{ 2  {\rm diag} (T {\rm diag} (c^*_s)c_s) -
                                    {\rm diag} (c^*_s)T {\rm diag} (c_s) \right\} -
\epsilon_i I,
\end{eqnarray*}

\begin{eqnarray*}
d^2 E_{ij} =
\frac{4}{h} \left\{ 4  {\rm diag} (c^*_i)T {\rm diag} (c_j) -
{\rm diag} (c^*_j)T {\rm diag} (c_i) -  {\rm diag} (T {\rm diag} (c^*_i)c_j) \right\}.
\end{eqnarray*}
%
%
Here, $d^2E_{ii}$ and $d^2E_{ij}$ are the diagonal and off-diagonal elements, respectively, of the block matrix of size $m \times m$.


As can be seen, multiplying a vector by such a matrix can be done in approximately $m^2 N \log_2 N$ operations. However, solving a system of linear equations with such a matrix is a rather complex task, as it requires an additional iterative process.

\subsubsection{Exact Transformations of the Discretized Hartree-Fock Energy Functional}

Let's consider the problem of finding the solutions to the Hartree-Fock equations for $m$ electron pairs in an orthonormal space of $N$ vectors $\{q_i\}$, where $N-m < m$. Such a problem may arise in relaxation-based minimization methods \cite{33}. In this case, the first $m$ vectors $\{q_i\}$ correspond to the solution of the Hartree-Fock equations at each iteration, while the remaining $N-m$ vectors represent directions for improving the solution. In this problem, there are $Nm$ unknowns and $\displaystyle \frac{m(m+1)}{2}$ constraints associated with the orthonormality of the solution. Therefore, if possible, we would like to simplify the problem.

It has been discovered that this problem can be transformed into a problem with $(N-m)m$ unknowns and correspondingly $\frac{\displaystyle (N-m)(N-m+1)}{\displaystyle 2}$ constraints. When $N-m << m$, this approach provides a significant advantage.

The most interesting cases are those where the values of $m$ range from 10 to 50. In such cases, a significant advantage can be observed when $N-m$ ranges from 1 to 10.

Note that in the case of a quadratic functional, this is evident, and the solutions are the $N-m$ largest eigenvalues and their corresponding eigenvectors. It turns out that a similar situation occurs in the case of the Hartree-Fock functional.

We will use the notations from equations (\ref{eq:3.4.2.1}) to (\ref{eq:3.4.2.2.2}). Let's find the vectors $c_i$, $i=1, \dots, m$ that correspond to the minimum of $E$. We construct an additional $N-m$ arbitrary vectors $\{c_i\}$ such that the resulting system of $N$ vectors remains orthonormal. Let $C=(c_1, \dots, c_N)$.

We transform the problem of minimizing $E$ with respect to the first $m$ vectors $c_i$ into the problem of minimizing the same $E$ with respect to the $N-m$ largest vectors $c_i$. Obviously, $CC^*=I$. Then the matrix $P=\{p_{ij}\}$ can be represented as $\displaystyle P = I - \sum^{N}_{i=m+1} c_i c^*_i$. Substituting this into (\ref{eq:3.4.2.2}):
\begin{eqnarray}
\nonumber
E = 2 \sum_{p,q=1}^N <p|H|q>
      \left( \delta_{pq} - \sum_{k=m+1}^N c_{pk}^* c_{qk} \right) +
\\
\nonumber
    \sum_{p,q=1}^N \sum_{r,s=1}^N G_{pqrs}
    \left(\delta_{pq} - \sum_{k=m+1}^N c_{pk}^* c_{qk}\right)
    \left(\delta_{rs} - \sum_{l=m+1}^N c_{rl}^* c_{sl}\right) =
\\
\label{eq:4.2.4.1}
\sum_{p=1}^N \left( 2 <p|H|q> +G_{pppp} \right) +
\\
\nonumber
\sum_{p,q=1}^N \sum_{k=m+1}^N \left\{ - 2 <p|H|q> +
    \sum_{r=1}^N (G_{pqrr}+G_{rrpq}) \right\} +
\\
\label{eq:4.2.4.2}
    \sum_{p,q=1}^N \sum_{r,s=1}^N G_{pqrs}
    \left( \sum_{k=m+1}^N c_{pk}^* c_{qk}\right) 
    \left( \sum_{l=m+1}^N c_{rl}^* c_{sl}\right).
\end{eqnarray}
%
%
As can be seen, $E$ consists of two terms: (\ref{eq:4.2.4.1}) and (\ref{eq:4.2.4.2}), where the term (\ref{eq:4.2.4.1}) is a constant and does not affect the value of $E$, while the term (\ref{eq:4.2.4.2}) is similar to the original problem but involves only $(N-m)m$ unknowns.



Let's consider another auxiliary problem: in the domain $D \in \R^3$, we have a three-dimensional grid with a step size of $h$ and size $n_i \times n_j \times n_k$. On this grid, we have piecewise polynomial $p$-smooth normalized functions $\{ S_{ijk} \}^{n_i, n_j, n_k}_{i,j,k=1}$ with finite support and periodic conditions along each Cartesian coordinate. Each of these functions is a tensor product of polynomials of degree up to $p+1$.

Let
$$
\phi=\sum_{i=1}^{n_i} \sum_{j=1}^{n_j} \sum_{k=1}^{n_k} a_{ijk} S_{ijk},
\hspace*{10mm}
\psi=\sum_{i=1}^{n_i} \sum_{j=1}^{n_j} \sum_{k=1}^{n_k} b_{ijk} S_{ijk}.
$$
%
%
It is evident that $\psi \phi \in \R^p$. Let's construct
$$ \xi=\sum_{i=1}^{n_i} \sum_{j=1}^{n_j} \sum_{k=1}^{n_k} \alpha_{ijk} S_{ijk} $$
%
%
such that the values of $\psi \phi$ and $\xi$ coincide at the grid nodes. Since $\xi \in \R^p$, the values of $\psi \phi - \xi$ and all its $p$ partial derivatives are zero at the grid nodes. Since the degree of the polynomials in $\psi \phi - \xi$ on each cube is at most $2p+2$, imposing $2p+2$ constraints allows us to express $\psi \phi - \xi$ as
$$ \psi \phi - \xi = \sum_{i=1}^{n_i} \sum_{j=1}^{n_j} \sum_{k=1}^{n_k}
  \beta_{ijk} {\cal S}_{ijk}, $$
%
%
where ${\cal S}_{ijk}$ is only defined in the cube $(i,j,k)$ and is identically zero in all other cubes. It is a tensor product of three one-dimensional functions $u_i(x)$, $u_j(y)$, and $u_k(z)$, each of which is a polynomial of degree $2p+2$ and has zero values and all $p-1$ derivatives on the cube's boundary.


This transformation allows us to represent any $\psi$ and $\phi$ as a sum with weights of $n_i n_j n_k$ pairs of regular basis functions, and it is an exact transformation. Using this transformation, we can modify equation (\ref{eq:4.2.1.2}) for any $p$ such that the block size of matrix $T$ becomes only $8 \times 8$. In this case, the computational cost of evaluating $V_{ee}$ will be only $32 m^2 N (\log_2 N + 16)$ operations.

\subsubsection{Approximate Computation of the Discretized Hartree-Fock Energy Functional}


Generalizing the previous result, we can:
\begin{itemize}
\item Neglect the term $\psi \phi - \xi$, then the computational cost of evaluating $V_{ee}$ will be only $\displaystyle \frac{1}{2} m^2 N \log_2 N$ operations.
\item Approximate $\psi \phi$ on a finer grid with step size $\displaystyle \frac{h}{K}$, where $K=2, 3, \dots$. In this case, the computational cost of evaluating $V_{ee}$ will be of the order $\displaystyle \frac{1}{2} m^2 N K \log_2 (NK)$ operations. The formula for estimating the accuracy of such an approximation can be derived for specific types of basis functions, but in general, we have not obtained a simple and clear expression.
\end{itemize}


Another simplification is related again to the term $V_{ee}$. Since the computational cost of evaluating $V_{en}$ and $T_e$ is only ${\tt const} N \log_2 N$ operations compared to ${\tt const} m^2 N \log_2 N$ for $V_{ee}$, this simplification is highly relevant. The idea is as follows: approximate $T$ from (\ref{eq:4.2.1.2}) as $\alpha I + \beta uu^T$, where $u=(1, \dots, 1)^T \in \C^N$. Then, due to the orthonormality of the solution, equation (\ref{eq:4.2.1.2}) can be rewritten as:

\begin{eqnarray}
\label{eq:4.2.5.1}
V_{ee} = \frac{1}{h} (\alpha z^* z + \beta (2 m - 1) m),
\end{eqnarray}
where $\displaystyle z = \sum_{s=1}^m H(c_s) c_s.$


Similarly, equation (\ref{eq:4.2.2.1}) can be rewritten as:
$$ dE_i =
\frac{2}{h^2} A c_i - \frac{4}{h} B c_i +
\frac{4}{h} \alpha  {\rm diag} (z) c_i + \frac{4}{h} \beta (2m - 1) -
\epsilon_i c_i.
$$
%
%
As a result, the computational cost of evaluating both the discretized Hartree-Fock functional and all its first partial derivatives reduces to ${\tt const} N \log_2 N$ operations. Since $T$ is a three-level block circulant symmetric matrix, this approximation can be performed as follows: $\displaystyle \min ||T - \alpha I - \beta u u^T ||_*$, in the Frobenius and spectral norms, requiring no more than $N \log_2 N$ arithmetic operations.


Based on experimental computations, the ratio
$$\frac{||T - \alpha I - \beta u u^T ||_*}{||T||_*}$$
%
%
varies within the range of $0.001$ to $0.01$, depending on the norm and the type of basis functions. This allows us to use such an approximation for rough calculations, initial approximations, or for efficient preconditioning.

Let's examine the dependencies
$$
\frac{||T_e||_F^2}{||V_{en}||_F^2} \hspace*{5mm} and \hspace*{5mm}
\frac{||T_e||_F^2}{||V_{ee}||_F^2}.
$$
%
%
Since $A$ is a three-level band circulant matrix, we have $||A||_F^2 = {\tt const} N$, where ${\tt const}$ is the sum of squares of elements in any row of matrix $A$. By construction, this quantity is of the order of unity. Therefore,
$$||T_e||_F^2 \simeq \frac{n^3}{h^4}.$$
From (\ref{eq:4.2.1.1}),
$$||B||_F^2 \le \sum_{s=1}^m Q_s^2 ||W||_F^2,$$
%
%
and, by construction, $||W||_F^2 \simeq n^2$. Considering that $n$ is inversely proportional to $h$, we obtain
\begin{eqnarray}
\label{eq:5.1}
\frac{||T_e||_F^2}{||V_{en}||_F^2} \ge C n^3 = C_1 N,
\end{eqnarray}
%
%
where $N$ is the size of matrices $T_e$ and $V_{en}$. It follows that as the number of unknowns in the discrete problem increases, the kinetic energy operator becomes larger in magnitude compared to $V_{ee}$.


A similar estimation can be obtained for $\displaystyle \frac{||T_e||_F^2}{||V_{ee}||_F^2}$. Since $||V_{ee}||_F^2 \le ||T||_F^2$, it follows that the kinetic energy operator becomes larger in magnitude compared to the nonlinear component of the Hartree-Fock equations as the size of the matrices increases, that is,

\begin{eqnarray}
\label{eq:5.2}
\frac{||T_e||_F^2}{||V_{ee}||_F^2} \ge C n^3 = C_2 N.
\end{eqnarray}
%
%
Similarly, estimates for these ratios in the spectral norm can be obtained. Then,
\begin{eqnarray}
\label{eq:5.3}
\frac{||T_e||_2}{||V_{en}||_2} \ge C_3 n
\hspace*{5mm} and \hspace*{5mm}
\frac{||T_e||_2}{||V_{ee}||_2} \ge C_4 n.
\end{eqnarray}
%
%
It can be easily observed that these estimates are optimal and cannot be improved.

\subsection{Solving Hartree-Fock Discrete Equations with Iterative Methods}


To find the minimum energy of Hartree-Fock, we applied the standard Kantorovich minimization procedure \cite{19,32}, commonly referred to as self-consistent field (SCF) method in the chemical literature. At each step, a linear eigenvalue problem for the Hartree-Fock operator is solved using the approximation from the previous step to form $V_{ee}$. The initial approximation is taken as the minimal eigenvalues and their corresponding eigenvectors of the linear part of the Hartree-Fock operator:
$$[T_e + V_{en}]\bar \phi_0 = \lambda \bar \phi_0,$$
$$[T_e + V_{en} + V_{ee}(\bar \phi_{i-1})]\bar \phi_i = \lambda \bar \phi_i, ~~ i=1, \dots.$$ 
%
%
Here, $\bar \phi_i$ is a set of refined eigenvectors corresponding to the minimal eigenvalues.


The convergence of this method for a single electron pair has been proven in \cite{36}, and it is conjectured that convergence also holds for multiple electron pairs.


Currently, the most popular methods for finding one or several minimal eigenvalues and their corresponding eigenvectors of a Hermitian matrix are the Lanczos methods with Chebyshev acceleration in the general case \cite{37}, and the Davidson method \cite{Dav} if a good preconditioner is known. From (\ref{eq:5.1}--\ref{eq:5.3}), the operator $T_e^{-1}$ serves as a good preconditioner for the Davidson method, as $||T_e||F^2$ grows linearly compared to $||V_{en}+V_{ee}||_F^2$ with respect to the problem size. Moreover, since $T_e^{-1}$ is also a three-level circulant matrix, the cost of constructing and multiplying with such a matrix is of the order $N \log_2 N$, which is significantly smaller ($m^2 N \log_2 N$) than multiplying with the entire Hartree-Fock operator. Furthermore, a symmetric three-level circulant matrix $C$ was used as a preconditioner, which was obtained from
$$\min_C || T_e + V_{en} + \frac{4}{h} \alpha  {\rm diag}  (z) + \frac{4}{h} \beta (2 m - 1) I - C ||_F.$$
%
%
It is evident that constructing $C$ is also possible with approximately $N \log_2 N$ arithmetic operations.


It can be shown that by shifting $C$ by an amount much smaller than $\frac{4}{h} \beta (2 m - 1)$, a matrix can be obtained with a condition number in the $||.||_2$ norm of at most $N$.


As it turns out, the ratio:
\begin{eqnarray}
\label{eq:prec}
\frac{|| T_e + V_{en} + \frac{4}{h} \alpha  {\rm diag}  (z) +s
                        \frac{4}{h} \beta (2 m - 1) I - C ||_F}{||C||_F}
\end{eqnarray}
%
%
is much smaller than unity and has a magnitude of $\sqrt{N}$, which means that on a finer grid, this ratio is even smaller. It has been shown in \cite{Dav} that the convergence of the Davidson algorithm depends on this quantity.




Based on these considerations, the Davidson method was chosen for finding the minimal eigenvalues and corresponding eigenvectors of the Hermitian matrix that arises at each step of the Kantorovich minimization procedure.

Currently, there are various variants of this method, so for clarity, we will present one of them and then explain why this particular method was chosen.

Let $A \in \mathbb{C}^{N \times N}$ be a Hermitian matrix for which we need to find the minimal eigenpair $(\nu, \psi)$. Let $P \simeq A$ be a preconditioner such that $(P-\gamma I)$ is easily invertible. Let $\psi_0$ be an initial approximation to such an eigenvector. Then, the iterations proceed as follows:

\bigskip

\noindent
$\nu_0=\psi_0^* A \psi_0$,

\noindent
$q_0 = \psi_0$,

\noindent
$\forall i=1, \dots, {\tt MAXSPACE}$:

$r_i = A \psi_{i-1} - \nu_{i-1} \psi_{i-1}$ --- computing the residual of the $i$-th iteration,

if $||r_i||_F < {\tt eps}$, a solution is found, and the iterations are stopped. Otherwise:

$q_i = (P - I\nu_{i-1})^{-1} r_i$,

orthogonalize $q_i$ with respect to all previous vectors $q_0, \dots, q_{i-1}$,

$H = (q_0, \dots, q_i)^* A (q_0, \dots, q_i)$,

find the minimal eigenvalue $\nu_i$ and its corresponding eigenvector $\alpha_i$ of the matrix $H$,

$\psi_i = (q_0, \dots, q_i) \alpha_i$

If $i={\tt MAXSPACE}$, set $\psi_0 = \psi_i$ and return to step 1.

\bigskip


During each iteration, it is necessary to store $i+1$ vectors from $\mathbb{C}^N$, perform approximately $4iN$ arithmetic operations during orthogonalization and the construction of the next approximation, and multiply a vector by the matrix twice.

An alternative to this method could be the Davidson method with a single matrix-vector multiplication per iteration, but it would require storing an additional $i$ vectors of the form $A (q_0, \dots, q_i)$ and performing $iN$ more arithmetic operations during orthogonalization.

If an approximate problem is being solved, where $A = T_e + V_{en}$, the arithmetic cost of multiplying a vector by a matrix ranges from $3.5N$ to $13.5N$. Therefore, using an algorithm with a single matrix-vector multiplication per iteration is not justified in terms of memory usage or the number of arithmetic operations performed.

However, if the vector is multiplied by the original Hartree-Fock operator without any simplifications, the number of arithmetic operations for each such multiplication can reach $500N$. In this case, the main computational bottleneck becomes the matrix-vector multiplication.

Thus, a dilemma arises: one algorithm requires half the memory, while the other works twice as fast. The slower algorithm was chosen because, for most problems, the data volume during iterations was the most critical resource.

Currently, there are two approaches for finding multiple eigenvectors. The first approach is to use the same Davidson method, but instead of seeking a single desired vector, multiple vectors are sought simultaneously. The second approach exploits the fact that at each iteration, the solution is orthogonalized to all the previously found eigenvectors.

In our case, the first approach was not feasible because the number of iterations in such a block method remained approximately the same, which would have resulted in a several-fold increase in memory usage. Therefore, only the second approach was used: orthogonalizing the next approximation to all the already found eigenvectors.

It turned out that the Davidson algorithm was suitable for this problem, as the number of iterations for test problems did not exceed 50.

Another simplification was made regarding the choice of the initial approximation. Since $V_{en}$ is reasonably well approximated by formula (\ref{eq:4.2.5.1}), an initial simplified nonlinear Hartree-Fock eigenvalue problem was solved using this approximate operator. The obtained approximation was then used to solve the general Hartree-Fock problem. In this case, the number of iterations for solving the general Hartree-Fock problem was reduced by 2-3 times.

\section{Approximation of the Solution and Hartree-Fock Operator}




Although the use of such generalized bases allows achieving any predefined accuracy, it requires significant computational costs.

At the same time, the solution of the Hartree-Fock equations often approximates a sum of single-atom solutions, centered and differently oriented in space.

Any such single-atom solution of the Schrödinger equation can be well approximated by the following expression:
\begin{equation}
\label{eq:6.1}
\psi(x,y,z) \simeq \sum_{r=1}^{R} \alpha_r \psi_r^{(x)}(x) \psi_r^{(y)}(y) \psi_r^{(z)}(z),
\end{equation}
%
%
where $R$ is a very small value. It will be referred to as the rank or tensor rank of such decomposition.




Based on numerous numerical experiments, it has been observed that the solutions of the Hartree-Fock equations also have a low-rank approximation. In other words, storing and utilizing each eigenvector of the Hartree-Fock equations does not require a three-dimensional object, which, when discretized as $200 \times 200 \times 200$, would require about 100 megabytes of memory for each eigenfunction. It is sufficient to have several one-dimensional objects of this kind without losing the accuracy of the approximation.

In fact, even with a discretization of the order of $1000 \times 1000 \times 1000$, it is enough to store only a few one-dimensional discrete vectors and their corresponding weight coefficients to obtain a good approximation of the wave function, which amounts to even less than one megabyte of data.

Let us reconsider Equation (\ref{eq:4.2.2.1}) and note that when using orthogonal bases, the vanishing of these first derivatives can be expressed as:

$$ F = \frac{2}{h^2} A - \frac{4}{h} B
+ \frac{8}{h} \sum_{s=1}^m {\rm diag}\left\{ T {\rm diag} (c_s^*) c_s \right\}
- \frac{4}{h} \sum_{s=1}^m \left\{ {\rm diag}(c_s) T {\rm diag}(c_s) \right\}
$$
$$ \forall i=1, ..., m: ~~ F c_i = \epsilon_i c_i.$$
%
%
Here, all the matrices used, such as $A$, $B$, and $T$, although being special (band or Toeplitz matrices), can also be represented as a sum of a finite and very limited number of tensor products of one-dimensional matrices. For example, $A$ is, by construction, a sum of three one-dimensional tensor matrices, $T$ can be approximated with an accuracy of $10^{-6}$ using only six one-dimensional tensor matrices, and depending on the arrangement of atoms, $B$ also has a very low tensor rank.






Multiplying a vector by such a tensor matrix results in a one-dimensional tensor. Moreover, for such multiplication, if the original matrix had dimensions of $\mathbb{C}^{n^3 \times n^3}$, we would only need to perform $\mathcal{O}(6n^2)$ arithmetic operations.

The use of Toeplitz or circulant structures slightly accelerates these multiplications. However, it should be noted that for most real problems, $n$ ranges from 100 to 1000, and the use of circulant matrices can ideally provide a speedup of 3-10 times, while for Toeplitz matrices, it is only a few times.

However, the use of arbitrary matrices, not necessarily Toeplitz or circulant, simplifies the problem of atom centering in the finite element.

When solving the Hartree-Fock problem for a predefined set of atoms with their given coordinates, it is generally not possible to construct a discretization grid in which each atom is precisely located at the center of a finite element, or it may require creating a grid with very small cells.

The use of tensor structure solves this problem. We can construct a tensor grid with a slightly non-uniform mesh such that each atom falls at the center of a finite element. An illustrative example of such a grid is shown in the following figure.

$$\includegraphics[scale=0.6]{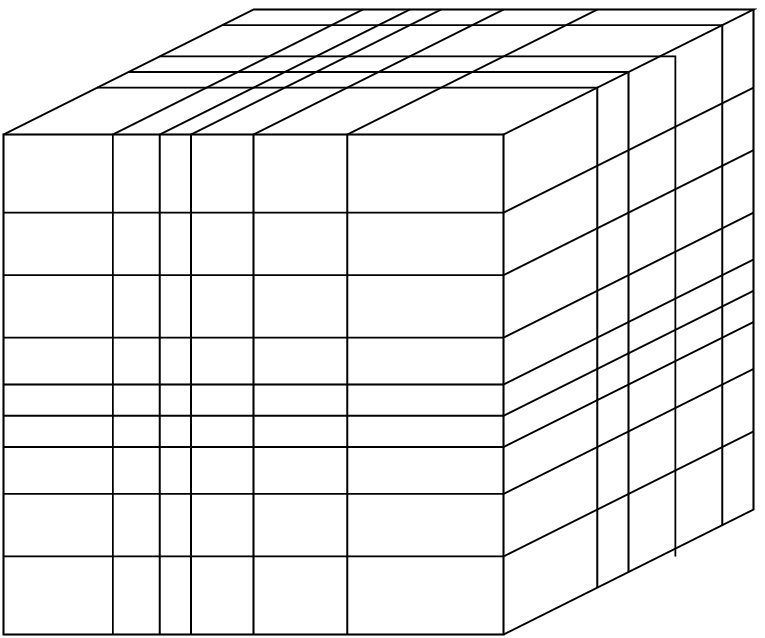}$$
\centerline{An example of the rectangular, non-uniform grid.}
\bigskip

It can be easily shown that with a given discretization step, in the worst case, we only need to double the step, which guarantees that each atom is always located at the center of a finite element, and the number of finite elements in each direction increases by no more than a factor of two.

When multiplying a sum of tensor matrices by a sum of tensor vectors, the tensor rank generally increases. However, since we expect the solution to have a low rank, it is necessary to compress the resulting tensor structures with a large rank to structures with a small tensor rank. Algorithms for efficiently performing such compression have been proposed by the author in \cite{3}.

Thus, if Toeplitz and circulant structures are used, the arithmetic complexity of each iteration of the Davidson method on uniform grids requires only $\mathcal{O}(n)$ arithmetic operations. If arbitrary irregular tensor grids are used, the complexity becomes $\mathcal{O}(n^2)$ arithmetic operations.

Of course, the constant factor in this asymptotic estimate is quite large and implicitly depends quadrati\-cally on the number of electrons. However, it is worth noting that for $N = 10^7$, the proposed method provides significantly higher accuracy than the classical Hartree-Fock method with an arbitrary basis based on one-atom molecular orbitals for small molecules with linear dimensions of no more than 10 Angstroms. In this case, $n \approx 200$, and the proposed method allows the use of only 100 MB to 1 GB of memory for the entire solution for molecular structures with tens of atoms and hundreds of electron pairs, along with a very small number of arithmetic operations, while achieving the desired accuracy of the solution.

In addition, in this work, the method was tested on a set of spatial structure calculations for several molecules from the HugeMDB database \cite{HugeMDB}, and the computational speed of solving the Hartree-Fock equations without and with the electron-electron interaction term, as well as the DFT method, was evaluated. The observed values on modern computers in 2023 are presented below.

Sets of molecules were used as test cases, each consisting of 100 random conformers with specified linear characteristics (L in Angstroms) and a given number of electron pairs (m). The calculations and properties of these calculations were performed for an 16-core AMD processor (approximately 640 Gflops/s peak performance in double precision) and are presented in the table below. For each set, the average values of rank (R), computation speed (T), and memory required (M) were calculated. In the table, these values are represented using a slash symbol, where the first value in each column corresponds to the characteristics for solving the Hartree-Fock equations without the electron-electron interaction term, the second value corresponds to the characteristics for solving the Hartree-Fock equations with the electron-electron interaction term, and the third value corresponds to the characteristics for solving the DFT equations.

$$\begin{tabular}{|c|c|c|c|c|} \hline
 L &   m &        R &              T & M     \\ \hline
 6 &  15 & 16/21/19 & 1.4m/5.3m/2.4m &  99MB/219MB/120MB \\ \hline
 8 &  25 & 18/24/20 &  3.9m/31m/8.7m & 216MB/518MB/248MB \\ \hline
10 &  40 & 18/25/21 &     9m/86m/17m & 405MB/1.1GB/510MB \\ \hline
12 &  60 & 19/25/22 &   22m/5.4h/45m & 0.7GB/1.8GB/0.8GB \\ \hline
15 & 100 & 18/24/22 &  1.5h/38h/3.8h & 1.3GB/3.5GB/1.4GB \\ \hline
20 & 160 & 18/23/20 &   6.5h/11d/20h & 2.5GB/7.0GB/2.7GB \\ \hline
\end{tabular}$$



\vspace*{-1mm}

\begin{thebibliography}{99}\setlength{\itemsep}{-0.36\baselineskip}

\bibitem{1} I. Ibragimov, Application of Structured Matrices for Solving Hartree-Fock Equations / Mat. Meth. and Comp. Moscow, INM RAS,
ISBN 5-201-08803-1, (1999): 144--173.

\bibitem{2} I.~Ibragimov, Algebraic structures and parallel computing in quantum chemical problems. PhD thesis, Moscow, INM RAS, (1999): 114.

\bibitem{19} W.~J.~Herbe, L.~Radom, P.~v.R.~Schleyer, J.~A.~Pople. {\it Ab initio} Molecular Orbital Theory. Springer-Verlag (1986): 327.

\bibitem{11} A.~P.~Seitsonen, M.~J.~Puska, and R.~M.~Neiminen. Real-space electronic-structure calculations: Combination of the finite-difference and conjugate gradient methods. Phys. Rev. B.: 51:20, (1995): 14057--14061.

\bibitem{12} M.~C.~Payne, M.~P.~Teter, D.~C.~Allan, T.~A.~Arias, and J.~D.~Joannopoulos. Iterative minimization techniques for {\it ab initio} total-energy calculations: molecular dynamics and conjugate gradients. Rev. of Mod. Phys. 64:4, (1992): 1045--1097.

\bibitem{20} A.~S.~Umar, M.~R.~Strayer, J.-S.~Wu, D.~J.~Dean, and M.~C.~G\"ucl\"u. Nuclear Hartree-Fock calculations with splines. Phys. Rev. C.: 44:6, (1991): 2512--2521.

\bibitem{24} J.-L.~Calais. Wavelets --- Something for Quantum Chemistry. Int. J. Quant. Chem. 58, (1996): 541--548.

\bibitem{26} G.~G.~Hall, and D.~Rees. A Discrete Look at Localization. Int. J. Quant. Chem. 53, (1995): 189--205.

\bibitem{34} I.V.~Ibraghimov, E.E.~Tyrtyshnikov. / Pro\-cee\-dings of the 9th International Conference ``Com\-pu\-ta\-tional Mo\-del\-ing and Computing in Physics''. 16-21.09.1996, D5:11-97-112, Dubna, (1997): 327--331.

\bibitem{31} C. de Boor. A Practical Guide to Splines. Springer New York, NY, 1978, 378p.

\bibitem{38} A.A.~Lokshin, S.L.~Lopatnikov, Yu.I.~Tarasov. The contraction mapping method in the symmetric eigenvalue problem. Moscow, MSU, (1995): 143.

\bibitem{32} L.~Heigeman, D.~Young. Applied iterative methods. Academic Press, (1981): 404.

\bibitem{33} V.~F.~Demianov, V.~N.~Malozemov. Introduction to minimax. Moscow, Nauka, (1972): 368.

\bibitem{36} A.A.~Lokshin, A.S.~Saakyan, Yu.I.~Tarasov. Parametric dependencies of eigenvalues. Moscow, MSU, (1997): 64.

\bibitem{37} B.~Parlett. The Symmetric eigenvalue problem. SIAM, (1998): 394.

\bibitem{Dav} M.~Crouzeix, B.~Philippe, M.~Sadkane. The Davidson Method. SIAM J. Sci. Comp., 15:1, (1994): 62--76. 

\bibitem{3} I. Ibragimov. A new approach to solving the problem of generalized singular value decomposition / Mat. Meth. and Comp.. Moscow, INM RAS, ISBN 5-201-08803-1, (1999): 193--201.

\bibitem{HugeMDB} HugeMDB: one of the largerst database of small molecules: https://www.multi-d.com/

\end{thebibliography}
\end{document}